\documentclass[12pt]{amsart}

\usepackage{amsmath, amssymb, graphics, setspace}
\usepackage[usenames]{color}
\input xypic

\usepackage{algorithm}
\usepackage[noend]{algpseudocode}

\algblock{Input}{EndInput}
\algnotext{EndInput}
\algblock{Output}{EndOutput}
\algnotext{EndOutput}

\makeatletter
\def\BState{\State\hskip-\ALG@thistlm}
\makeatother

\usepackage{verbatim}
\usepackage{url}
\usepackage{bm}
\usepackage{tikz}
\usepackage{appendix}
\usepackage[noend]{algpseudocode}

\newcommand{\mathsym}[1]{{}}
\newcommand{\unicode}[1]{{}}

\makeatletter
\def\BState{\State\hskip-\ALG@thistlm}
\makeatother

\usepackage{amsmath}
\usepackage[centertags]{amsmath}
\usepackage{amsfonts}
\usepackage{amssymb}
\usepackage{amsthm}
\usepackage{newlfont}
\usepackage{url}
\usepackage{bm}
\usepackage{chngcntr}
\usepackage{mathtools}
\usepackage{blkarray, bigstrut} %

\theoremstyle{plain}
\newtheorem{thm}{Theorem}

\newtheorem{lem}[thm]{Lemma}
\newtheorem{lem*}[thm]{Lemma}
\newtheorem{prop}[thm]{Proposition}

\theoremstyle{definition}
\newtheorem{dfn}{Definition}
\theoremstyle{remark}
\newtheorem{rem}{Remark}
\newtheorem{rem*}{Remark}
\newtheorem{ex}[rem]{Example}
\newtheorem{example}[rem]{Example}

\numberwithin{rem}{section} %{subsection}
\numberwithin{dfn}{section} %{subsection}
\numberwithin{equation}{section} %{subsection}
\numberwithin{thm}{section} %{subsection}

\def\!{\operatorname{!}}

\def\C{\mathbb C}
\def\F{\mathbb F}
\def\G{\mathbb G}

\def\R{\mathbb R}

\def\1{\bold 1}

\def\diag{\operatorname{diag}}
\def\deg{\operatorname{deg}}

\def\Spec{\operatorname{Spec}}

\def\Hom{\operatorname{Hom}}
\def\GL{\operatorname{GL}}

\def\Ext{\operatorname{Ext}}

\def\rk{\operatorname{rk}}

\def\deg{\operatorname{deg}}

\def\End{\operatorname{End}}

\usepackage{amssymb}

\usepackage{enumitem}

\usepackage{graphicx}
\usepackage{xy}
\usepackage{float}
\usepackage[T1]{fontenc}

\makeatletter
\@namedef{subjclassname@2020}{%
	\textup{2020} Mathematics Subject Classification}
\makeatother

\theoremstyle{definition}

\numberwithin{equation}{section}

\newcommand{\Der}{\mathrm{Der} }

\newcommand{\Derin}{\Der_{in}}

\newcommand{\lra}{\longrightarrow}

\newcommand{\podwzorem}[2]{\underbrace{#1}\limits_{#2}}

\begin{document}

	\baselineskip=17pt
	
	%%%%%%%%%%%%%%%%
	
	\title{Weil-Barsotti formula for  $\mathbf{T}$-modules}

	\author{Dawid E. K{\k e}dzierski, Piotr Kraso{\'n}$^\star$}
	
	\thanks{$^\star$Corresponding author.}
	\date{\today}
	
	\address{  Institute of Mathematics, Department of Exact and Natural Sciences, University of Szczecin, ul. Wielkopolska 15, 70-451 Szczecin, Poland 
	}
	\email{dawid.kedzierski@usz.edu.pl}
	\email{piotrkras26@gmail.com}

	
\thanks{}

\maketitle

\newcommand{\tm}{$\mathbf{t}$-}
\newcommand{\tsm}{$\mathbf{t}^{\sigma}$-}
\newcommand{\functor}{\mathrm{Rep}}

%\begin{document}
\begin{abstract}
In the work of M. A. Papanikolas and N. Ramachandran [A Weil-Barsotti  formula for Drinfeld modules, Journal of Number Theory 98, (2003), 407-431] the Weil-Barsotti formula for the function field case  concerning $\Ext_{\tau}^1(E,C)$ where $E$ is a  Drinfeld module  and $C$ is the Carlitz module was proved. We generalize this formula  to the case where $E$ 
is a strictly pure \tm module $\Phi$ with the zero nilpotent matrix $N_\Phi.$ For such a \tm module $\Phi$ we explicitly compute its dual \tm module ${\Phi}^{\vee}$ as well as its double dual  ${\Phi}^{{\vee}{\vee}}.$  This computation is done in a  subtle way   by  combination of  the \tm reduction algorithm developed by  F. G{\l}och, D.E. K{\k e}dzierski, P. Kraso{\'n} [ Algorithms for determination of \tm module structures on some  extension  groups , 
arXiv:2408.08207] and the methods of the work of D.E. K{\k e}dzierski and  P. Kraso{\'n} [On $\Ext^1$ for Drinfeld modules, Journal of Number Theory 256 (2024) 97-135].  We  also 
give a counterexample to the Weil-Barsotti formula if the nilpotent matrix $N_{\Phi}$ is non-zero. 	\end{abstract}

\section{Introduction}
Duality is one of the fundamental concepts  that play an important role  in various branches of modern mathematics.  The idea of (concrete) duality  relies on associating  an object X with another   object $C(X)$  in the same category ${\cal C}$ which is the space of morphisms from $X$ to a simpler object $C$ i.e. $C(X)={\Hom}_{\cal C}(X,C).$
This kind of idea was explored by I. Gelfand \cite{ge1,ge2} who defined a map $A\rightarrow C(\Spec A)$ which connected  the space of  continuous (in a suitable topology) functions on the space of multiplicative linear forms on $A$  to a Banach algebra $A.$ The space of such forms  may be  naturally identified with the space of maximal ideals of $A.$
Then  A. Grothendieck (cf. \cite[p.397-398]{ca1} transferred this idea to the algebraic situation, associating  a commutative ring $A$ with
the space  of sections of the structure sheaf ${\Gamma}(\Spec A,\cal O).$ For the evolution of the various concepts  of duality and examples of these, the reader can consult \cite{kc14}.

Let $\cal A/F$ be an abelian variety over a field $F$ and $\G_m/F$ be the multiplicative algebraic group. Let ${\cal A}^{\prime}$ be the dual to an  abelian variety $\cal A$. 
The classical Weil-Barsotti formula asserts the  following natural functorial isomorphism \cite[Chapter VII, sec. 16, Theorem 6]{se1}, \cite[III. 18 ]{oort}:
\begin{equation*}
{\cal A}^{\prime}(F)\cong\Ext_F^1({\cal A},\G_m)
\end{equation*}
where $\Ext_F^1({\cal A},\G_m)$ denotes the group  of  extensions of algebraic groups over $F$:
\begin{equation*}
0\rightarrow \G_m \rightarrow E \rightarrow  \cal A \rightarrow 0
\end{equation*}
with the Baer addition (cf. \cite{hs} ).  The Cartier-Nishi theorem states that there is a canonical isomorphism \cite{ca},\cite{nishi}:
\begin{equation*}
{\cal A}(F)\cong \Ext_F^1({\cal A}^{\prime},G_m).
\end{equation*}
This can be restated as $({\cal A}^{\prime})^{\prime}\cong A.$

There is a deep analogy between the theory of elliptic curves over complex numbers and the theory of Drinfeld modules (cf. \cite[sec2.5]{bp20}) and more generally between the theories of abelian varieties and \tm modules \cite{g96}. 

In \cite{ta}, Taguchi  defined the analog of the Cartier dual $E^{\vee}$  for a finite \tm module $E$ and constructed the Galois compatible Weil pairing on the torsion points of  the Drinfeld modules $E$ and $E^{\vee}$  providing an example of such an analogy. He also remarks \cite[Remark 5.2]{ta} that this pairing can also  be  defined for \tm modules which we  call 
(cf. sec. 2) strictly pure with no nilpotence. It is well known that the category of \tm modules is anti-equivalent with the category of \tm motives \cite{a},\cite{gl3}. An attempt  to define a duality for  \tm motives corresponding to such \tm modules was  undertaken in \cite{gl}. Our strategy in this paper  is different. Our main goal is to generalize the 
main theorem of \cite{pr}. There in,
  M. Papanikolas and  N. Ramachandran    proved the Weil-Barsotti formula for Drinfeld modules as well as the analog of the Cartier-Nishi theorem. 
  We slightly generalize their definition of duality.
  Then  we have to prove  that the double  dual ${\Phi}^{\vee\vee}$  of a \tm module $\Phi$ exists and is isomorphic to ${\Phi}$ for the class of \tm modules under consideration. The main problem is the fact that the dual ${\Phi}^{\vee}$  of  a strictly 
  pure \tm module with no nilpotence is not necessarily strictly pure. 
 (All necessary definitions needed for understanding Theorem \ref{wb11} and Theorem \ref{main}  are given in section 2.)
One of the main  results from \cite{pr} reads as follows:
\begin{thm}[\cite{pr}]\label{wb11}
Let $E$ be a Drinfeld module of rank $r \geq 2.$
\begin{enumerate}
\item[(a)] The group $\Ext^1(E,C)$ is naturally a \tm module of dimension r and sits in an
exact sequence of t-modules
$$0\rightarrow E^{\vee}\rightarrow \Ext^1(E,C) \rightarrow \G_a\rightarrow 0.$$
Furthermore, $E^{\vee}$ is the Cartier-Taguchi dual \tm module associated to $E$ \cite{ta}, and in particular, $E^{\vee}$ is isomorphic to the $(r -1)$-st exterior power ${\wedge}^{r-1}E$ of $E.$
\item[(b)] The group $\Ext^1(E^{\vee},C)$ is also naturally a \tm module of dimension $r $ and sits in an exact sequence
$$0\rightarrow E\rightarrow \Ext^1(E^{\vee},C) \rightarrow \G_a^{r-1}\rightarrow 0.$$
Moreover, we have a biduality: $(E^{\vee})^{\vee}\cong E.$
\item[(c)] Any morphism $\beta : E\rightarrow F$ of Drinfeld modules (of rank $\geq 2$) induces a
morphism of dual \tm modules ${\beta}^{\vee}: F^{\vee} \rightarrow E^{\vee}.$
\end{enumerate}
\end{thm}
Here $C$ denotes the Carlitz module.
The above theorem was proved by identifying the space of extensions $ \Ext^1(E,C)$ as the quotient of the space of biderivations by the space of inner biderivations (see sec.2).
\begin{rem}\label{not}
Following the notation in \cite{kk04}, the extension group in the category of \tm modules will be denoted as ${\Ext}_{\tau}^1(\Phi,\Psi)$  rather than  ${\Ext}^1(\Phi,\Psi).$ This is to indicate that 
extensions are considered in the additive subcategory of \tm modules sitting inside the abelian category of ${\mathbb F}_q[t$]-modules (cf. also Remark \ref{sub}).
\end{rem}

In \cite{kk04} we gave a general procedure for establishing the \tm module structure on ${\Ext}_{\tau}^1(\Phi,\psi)$ where $\Phi$ is a {\it strictly pure} \tm module and $\psi$ is a Drinfeld module under the hypotheses that $\deg_\tau\Phi>\rk\psi.$
In \cite{gkk} we built a general algorithm, called the reduction algorithm, for determining the \tm module structure on ${\Ext}_{\tau}^1(\Phi,\Psi),$ where $\Phi$ is a strictly pure \tm module 
and $\deg_\tau\Phi>\deg_\tau\Psi.$

In this paper we compute the \tm module structure of ${\Phi}^{\vee}:=\Ext_{0,\tau}^1(\Phi,C)$ for $\Phi$ a strictly pure  \tm module with no nilpotence (i.e. $N_{\Phi}=0$) and $\deg_\tau\Phi\geq 2.$
 We also  explicitly  compute ${\Phi}^{{\vee}{\vee}}=\Ext_{0,\tau}^1(\Phi^{\vee},C).$ As $\Phi^{\vee}$ is usually not a strictly pure \tm module, special care is needed for performing  consecutive reductions. 
  In fact to reduce correctly  one has to carefully choose the order of  reductions. 
 Our analog of the Weil-Barsotti formula is the following theorem:
 \begin{thm}\label{main}
 Let $\Phi$  be  a strictly pure \tm module of dimensions $d$ with no nilpotence. If  $n=\deg_\tau\Phi\geq 2$ then  
		\begin{itemize}
			\item[$(a)$] $\Ext^1_{\tau}(\Phi,C)$ has a natural structure of a  \tm module of dimension $n\cdot d$
			and sits in the following exact sequence of \tm modules
			$$0\lra {\Phi}^{\vee}\lra \Ext^1_{\tau}(\Phi, C )\lra \G_a^{d}\lra0,$$
\item[(b)] The group $\Ext_{\tau}^1(\Phi^{\vee},C)$ is also naturally a \tm module of dimension $d\cdot n $ and sits in an exact sequence
$$0\rightarrow \Phi\rightarrow \Ext_{\tau}^1(\Phi^{\vee},C) \rightarrow \G_a^{(n-1)d}\rightarrow 0.$$
Moreover, we have a biduality: $({\Phi}^{\vee})^{\vee}\cong\Phi.$
\item[(c)] Any morphism $\beta : \Phi\rightarrow \Psi$ of strictly pure \tm modules (of degree $\geq 2$) with no nilpotence  induces a
morphism of dual \tm modules ${\beta}^{\vee}: \Psi^{\vee} \rightarrow \Phi^{\vee}.$
\end{itemize}
 \end{thm}
 
 \begin{rem}\label{special}
 Notice that Theorem \ref{wb11} is the specialization of Theorem \ref{main} for $d=1.$ For $d>1$ we do not   claim that ${\Phi}^{\vee}$ is the external power of ${\Phi}.$
 In  \cite[sec. 1.3]{gl1} the authors remark that the external power of a  \tm module (or rather a \tm motive) without nilpotence is almost always a \tm module with non-zero 
 nilpotent part. The only exception is for the case $d=1.$ But  we determined that for our class of \tm modules the dual \tm module has no nilpotence cf. (\ref{mac}) and thus cannot be an external power of the original \tm module if $d>1.$
 \end{rem}
 
As we mentioned above, section 2 contains the definitions and facts needed further in the paper. In section 3 we  perform necessary computations and prove Theorem \ref{main}. In section 4 we give a counterexample to the Weil-Barsotti formula where the nilpotent matrix $N_{\Phi}\neq 0.$ We also added an appendix in which we compare our dual with that of Y. Taguchi. We proved an analogous result concerning the existence of Weil pairing for 
a stricly pure \tm module with no nilpotence.

\section{basic definitions and facts}
In this section we give only basic facts concerning \tm modules. We limit ourselves only to the algebraic side of the theory of \tm modules. The reader is advised to consult  excellent sources like \cite{g96}, \cite{MP}, \cite{th04}. \cite{f13}.

Let $A={\mathbb F}_q[t]$ be a polynomial ring over a finite field ${\mathbb F}_q$ with $q$-elements 
and let $k={\mathbb F}_q(t)$ be the quotient field of $A.$
Let  $v_{\infty}: k\rightarrow {\mathbb R}\cup\{\infty\}$ be the normalized
valuation associated to $\frac{1}{t}$ (i.e. $v_{\infty}(\frac{1}{t})=1$). Let $K$  be a completion of $k$ with respect to $v_{\infty}$ and let $\overline K$ be its fixed  algebraic closure. It turns out that ${\overline K}$ is neither complete nor 
locally compact. 
By  ${\mathbb C}_{\infty}$  we denote  the completion of $\overline K$ with respect to the metric induced by   the extension $\overline{v}_{{\infty}}$  of  the valuation $v_{\infty}.$

An $A$-field $L$ is a fixed morphism  
 ${\iota}: A \rightarrow L$. The kernel of $\iota$ is a prime ideal $\cal P$ of $A$ called the characteristic. The characteristic of $\iota$ is called finite if ${\cal P}\neq 0,$ or generic (zero) if ${\cal P}= 0.$  We denote  $\theta:=\iota(t).$ 

   The endomorphism ring ${\End}(\G_{a,L})$ where $\G_{a,L}$ is the additive algebraic group over L, is the skew polynomial ring $L\{\tau\}.$  
  The endomorphism $\tau$ is 
the map $x\rightarrow x^q$ and therefore one has the commutation relation $\tau x=x^q\tau$ for $x\in L.$

The following definition was introduced by G. Anderson \cite{a}:
\begin{dfn}
A ${d}$-dimensional \tm module  over $L$ 
 is an algebraic group $E$ defined over $L$ and isomorphic over $L$ to $\G^{d}_{a,L}$  together with a choice of ${\mathbb{F}}_q$-linear endomorphism $l:E\rightarrow E$ such that $\partial(l-\theta)^n{\mathrm{Lie}}(E)=0$ for a  sufficiently large $n$. 
  Notice that 
$\partial(\cdot)$ stands here for the differential of an endomorphism of algebraic groups. The choice of  an isomorphism $E\cong \G^{d}_{a,L}$ after transferring the aforementioned endomorphism $l$  yields  the  
${\mathbb F}_q$ -algebra homomorphism
\begin{equation*}
{\Phi} : {\mathbb F}_q[t]\rightarrow {\mathrm{Mat}}_m (L\{\tau\})
\end{equation*}
such that ${\Phi}_t:={\Phi}(t)$, as a polynomial in $\tau$  with coefficients in ${\mathrm{Mat}}_m(L)$  is of the following form:
\begin{equation}\label{eq:def:t-modul}
{\Phi}_t=({\theta}I_d+N_{\Phi})\tau^0+M_1{\tau}^1+\dots +M_r\tau^r
\end{equation}
where $I_d $ is the identity matrix and $N_{\Phi}$ is a nilpotent matrix. 
\end{dfn}
We will denote a \tm module by $(E,\Phi)$ or simply by $\Phi$.
It is clear that a \tm module  $\Phi$  is uniquely determined by the value ${\Phi}_t.$
We  also  consider the  zero \tm module given by the map ${\mathbb F}_q[t]\rightarrow 0$. 
The category of \tm modules along with the zero \tm module becomes additive.

\begin{dfn} 
	Let $(E,\Phi)$ and $(F,\Psi)$  be  two \tm modules of dimension $d$ and $e$, respectively. A morphism $f:(E,\Phi)\lra (F,\Psi)$
	is then a homomorphism of algebraic groups over $L$ that preserves the chosen endomorphisms i.e.
	the following diagram commutes:
	\begin{equation}\label{diagram hom}
\xymatrix{
\G_{a,L}^d\ar[r]^f\ar[d]^{t_{\Phi}} & \G_{a,L}^e\ar[d]^{t_{\Psi}} \\
\G_{a,L}^d\ar[r]^f& \G_{a,L}^e
}
\end{equation} 
where for a \tm module $\Phi$ given by  \eqref{eq:def:t-modul}, $t_\Phi : \G_a^d \rightarrow \G_a^d$ is the following map of algebraic groups:
$$t_\Phi\big(  X \big)= ({\theta}I_d+N_{\Phi})X+M_1X^q+\dots +M_rX^{q^r},$$
where $X=[X_1,\dots X_d]^T\in \G_a^d.$
\end{dfn}	
Again after the choices of basis of $\G^d_{a,L}$ and $\G^e_{a,L}$ a morphism of \tm modules over $L$ is given by a matrix $f\in\mathrm{Mat}_{d\times e}(L\{\tau\})$
	 such that
	$$
	f\Psi_t = \Phi_tf.
	$$

A Drinfeld module is a \tm module of dimension 1
and is given  by  a homomorphism of ${\mathbb F}_{q}$-algebras $\phi : A\rightarrow L\{\tau\},$ $ a\rightarrow {\phi}_{a}\,,$  such that 
\begin{enumerate}
\item[1.] $\partial\circ {\phi}={\iota},$
\item[2.] for some $a\in A, \, {\phi}_{a}\neq {\iota}(a){\tau}^{0},$
\end{enumerate}
where $\partial\big({\sum}_{i=0}^{i={\nu}}\,\,a_{i}{\tau}^{i}\big)=a_{0}. $ The characteristic of a Drinfeld module is the characteristic of ${\iota}.$

Drinfeld modules were introduced by V. Drinfeld in \cite{d74} who called them {\em elliptic modules}. 
The simplest Drinfeld module appeared in  L. Carlitz's 1935 paper \cite{c35}.
This Drinfeld module $C:  A\rightarrow L\{\tau\}$  determined by $C_t=\theta +\tau$ is called the Carlitz module.

A \tm module $\Phi_t=({\theta}I_d+N_{\Phi})\tau^0+M_1{\tau}^1+\dots +M_r\tau^r$ is called strictly pure if $M_r$ is an invertible matrix (cf.\cite[subsection 1.1]{np}).
A \tm module  $\Phi_t=({\theta}I_d+N_{\Phi})\tau^0+M_1{\tau}^1+\dots +M_r\tau^r$ is said to have no nilpotence  if $N_{\Phi}=0$.

Over $\C_\infty$
  the rank of $\Phi$ is defined as the rank of the period lattice of $\Phi$ as a $\partial{\Phi}(A)$-module (cf. \cite[Section t-modules]{bp20}). 
  For a Drinfeld  module its rank equals $\deg_\tau\phi_t$. As this is not necessarily the case for  $d>1,$ we  additionally define the  $\tau$-degree  of a \tm module:  $\deg_\tau\Phi:=\deg_\tau\Phi_t$.   
    
\begin{rem}\label{sub}
Notice that  a \tm module $(E,\Phi)$ of  dimension $d$ induces an   $\F_q[t]-$module  structure on ${\tilde{L}}^d$ for any algebraic extension  ${\tilde{L}}$ of $L.$
$$a\cdot x = \Phi_a(x)\quad \textnormal{for}\quad a\in\F_q[t], \quad x\in {\tilde L}^d.$$
Such a module   is called the Mordell-Weil group ${\Phi}({\tilde L})$.  Moreover,  any morphism of \tm modules $f:(E,\Phi)\lra (F,\Psi)$ induces a morphism  of $F_q[t]-$ modules 
$f: \Phi({\tilde L})\lra \Psi({\tilde{L}})$. Then the category of \tm modules can be considered as an additive subcategory of the abelian category of ${\mathbb F}_q[t]-$modules. This 
subcategory is not full (cf. \cite[Example 10.2]{kk04}). Therefore the $\Hom$-set in the category of \tm modules will be denoted 
as ${\Hom}_{\tau}$ i.e. $\Hom(\Phi,\Psi):={\Hom}_{\tau}(\Phi, \Psi).$ 

Similarly, by  $\mathrm{Ext}^1_{\tau}(\Phi, \Psi)$ we denote the Bauer group of extensions of \tm modules i.e. the group of exact sequences
\begin{equation}\label{seq}
0\rightarrow (F,\Psi)\rightarrow (Y ,\Gamma) \rightarrow (E,\Phi) \rightarrow 0
\end{equation}
with the usual addition of extensions (cf. \cite{hs}). 
\end{rem}
An extension of a $\mathbf t$-module $\Phi:\F_q[t]\lra {\mathrm{Mat}}_d(L\{\tau\})$ by $\Psi:\F_q[t]\lra {\mathrm{Mat}}_e(L\{\tau\})$ can be determined by a biderivation, i.e. $\F_q-$linear map $\delta:\F_q[t]\lra {\mathrm{Mat}}_{e\times d}(L\{\tau\})$ such that
\begin{equation*}	
	\delta(ab)=\Psi(a)\delta(b)+\delta(a)\Phi{(b)}\quad \textnormal{for all}\quad a,b\in\F_q[t].
\end{equation*}
The $\F_q-$vector space of all biderivations will be denoted by $\Der(\Phi, \Psi)$ (cf. \cite{pr}). 
 The map $\delta\mapsto \delta(t)$ induces the $\F_q-$linear isomorphism of the vector spaces $\Der(\Phi, \Psi)$ and $ {\mathrm{Mat}}_{e\times d}(L\{\tau\})$.
  Let $\delta^{(-)}: {\mathrm{Mat}}_{e\times d}(L\{\tau\})\lra \Der(\Phi, \Psi)$ be an $\F_q-$linear map defined by the following formula:
  \begin{equation}\label{delta}
  	\delta^{(U)}(a)=U\Phi_a - \Psi_aU\quad \textnormal{for all}\quad a\in \F_q[t]\quad\textnormal{and}\quad U\in {\mathrm{Mat}}_{e\times d}(L\{\tau\}).
  \end{equation}
	The image of the map  $\delta^{(-)}$ is denoted by $\Derin(\Phi, \Psi)$, and is  called the space of inner biderivations. 
We have the following $\F_q[t]-$module isomorphism (cf.  \cite[Lemma 2.1]{pr}):
	\begin{align*}
		\Ext^1_{\tau}(\Phi,\Psi)\cong\mathrm{coker}\delta^{(-)}=\Der(\Phi, \Psi)/\Derin(\Phi, \Psi).
	\end{align*}
Recall that an $\F_q[t]-$structure on the quotient of the space of   biderivations is given by the following formula:
$$a*\Big(\delta+\Derin(\Phi,\Psi)\Big):= \Psi_a\delta+\Derin(\Phi,\Psi)$$
  Later on, we omit  $+\Derin(\Psi,\Phi)$ when we consider the coset  $\delta+\Derin(\Phi,\Psi)$.
	
	For $V\in {\mathrm{Mat}}_{n_1\times n_2}(L\{\tau\})$ let $\partial V\in {\mathrm{Mat}}_{n_1\times n_2}(L)$    be the constant term of  
$V$ viewed as a polynomial in $\tau.$ 
For \tm modules $\Phi$ and $\Psi$ let
\begin{equation*}
\Der_0(\Phi,\Psi)=\{ \delta \in \Der(\Phi,\Psi)\,\,\mid \,\, \partial\delta (t)=0\}
\end{equation*}
and
\begin{equation*}
	\Der_{0,in}(\Phi,\Psi)=\{ \delta \in \Derin(\Phi,\Psi)\,\,\mid \,\, \partial\delta (t)=0\}.
\end{equation*}
As in  \cite[p. 413]{pr} we make  the following definition:
\begin{dfn}\label{zero} 
For any \tm modules $\Phi$ and $\Psi$
\begin{equation*}
\Ext^1_{0,\tau}(\Phi, \Psi):=\Der_0(\Phi,\Psi)/ \Der_{0,in}(\Phi,\Psi). 
\end{equation*}
\end{dfn}
In \cite[Corollary 2.3]{pr} for any \tm modules ${\Phi}$ and $\Psi,$ the exactness of the following  sequence  of $F_q[t]$-modules: 
$$0\lra \Ext_{0,\tau}^1(\Phi,\Psi)\lra \Ext^1_{\tau}(\Phi,\Psi)\lra \Ext^1(Lie(\Phi), Lie(\Psi))$$
 was proved. Moreover if $\theta:=\iota(t)$ is transcedental over ${\mathbb F}_q$ then the last arrow was proved to be   an epimorphism  (see  loc. cit.).
 
 Sometimes the $\F_q[t]-$  module structure  on $\Ext_\tau^1(\Phi,\Psi),$ and therefore also on  $\Ext_0(\Phi,\Psi),$ comes from the \tm module action.
Two such cases were studied  in \cite{pr}. The first case is  when 
$\Phi$ is a Drinfeld module of rank greater or equal to two and $\Psi$ is the Carlitz module. The second is when     
 $\Phi$ and $\Psi$  are  tensors of Carlitz modules  under the assumption  that the dimension of  $\Phi$  is smaller than that of $\Psi$. 
 Some further examples were given in
 \cite{kk04}, and finally in  \cite{gkk} an algorithm (called the   \tm reduction algorithm) was constructed. This algorithm allows one to determine the \tm module structure on 
  $\Ext_\tau^1(\Phi,\Psi)$ for wide classes of  \tm modules $\Phi$ and $\Psi$. 

\begin{dfn}\label{dual2}
For a \tm module $\Phi$ we define its dual by the formula:
\begin{equation*}
{\Phi}^{\vee}=\Ext^1_{0,\tau}(\Phi, C).
\end{equation*}
\end{dfn}
\begin{rem}\label{remdd}
In general, if ${\Phi}$ is not a strictly pure \tm module then ${\Phi}^{\vee}$ is an ${\mathbb F}_q[t]$-module  which is not necessarily a \tm module (cf. Example \ref{exs}).
\end{rem}

 \section{Main results}\label{Main results}
 In this section we consider a strictly pure \tm module  $\Phi$ of dimension $d$ with no nilpotence given by
 \begin{equation}\label{eqq}
  \Phi_t=\theta I_d+\sum\limits_{i=1}^{n}A_i\tau^i,\quad \textnormal{where}\quad n>1.
  \end{equation}
   Thus $A_n$ is an invertible matrix.
 We adapt the following notations:
\begin{itemize}
\item[(i)]$E_{1\times j}$ be  the matrix of type $1\times d$ with the only nonzero entry, equal to one, at the place $(1,j)$.
\item[(ii)] $B_0=A_n^{-1},\,\, B_j=A_n^{-1}A_j,\,\, j=1,2,\dots ,n.$ Notice that $B_n$ is the identity matrix $I_d.$
\item[(iii)] $B_j=[b_{j,k\times l}]_{k,l},$ where $b_{j,k\times l}$ denotes the $(k, l)$
entry of $B_j$.

\item[(iv)] $L\{\tau\}_{\langle 1, n )}=\Big\{ \sum\limits_{i=1}^{n-1} c_{i}\tau^{i}\mid c_i\in L \Big\},$
\item[(v)] 
 $A_1^{\vee}(B_1,\dots , B_{n-1})= \begin{bmatrix}
C_{1\times 1}& C_{1\times 2}&\cdots  & C_{1\times d}\\
C_{2\times 1}& C_{2\times 2} &\cdots  & C_{2\times d}\\
\cdots& \cdots&\cdots  & \cdots\\
C_{d\times1}& C_{d\times 2} &\cdots  & C_{d\times d}\\
\end{bmatrix}$

where  $$C_{i\times i} =\begin{blockarray}{*{7}{c} l}
    \begin{block}{*{7}{>{$\footnotesize}c<{$}} l}
     &  &  & & & & \\
   \end{block}
   \begin{block}{[*{7}{c}]>{$\footnotesize}l<{$}}
    0 & 0& \cdot & \cdot &\cdot & 0&  -b_{1,i\times i} &         \\
          1 & 0 & 0 & \cdot & \cdot & \cdot & -b_{2,i\times i} &\\
               0 & 1 &  0 & \cdot & \cdot &\cdot & \cdot &\\
                 \cdot & \cdot & \cdot & \cdot &\cdot & \cdot & \cdot  &\\
                   \cdot & \cdot & \cdot & \cdot &\cdot & 0& \cdot  &\\
                    \cdot & \cdot & \cdot & \cdot &\cdot & 0 & -b_{n-2,i\times i}  &\\
                      \cdot & \cdot & \cdot & \cdot &0 & 1 & -b_{n-1,i\times i}  &\\
                        \end{block}
  \end{blockarray}$$ {and}  
  $$C_{i\times j}= \begin{blockarray}{*{7}{c} l}
    \begin{block}{*{7}{>{$\footnotesize}c<{$}} l}
     &  &  & & & & \\
    \end{block}
    \begin{block}{[*{7}{c}]>{$\footnotesize}l<{$}}
     0& \cdot & \cdot & \cdot &\cdot & 0 &-b_{1,j\times i} &\\
   \cdot& \cdot &   \cdot & \cdot&\cdot & \cdot& -b_{2,j\times i} &  \\
   \cdot & \cdot & \cdot & \cdot &\cdot & \cdot & \cdot   &\\
       \cdot & \cdot & \cdot &\cdot & \cdot & \cdot & \cdot &\\
      0 & \cdot & \cdot & \cdot & \cdot &0 &  -b_{n-1,j\times i} & \\
\end{block}
  \end{blockarray}
  \quad {\mathrm{for}}\quad i\neq j.
$$
\item[(vi)] $A_2^{\vee}(B_0)= \begin{bmatrix}
D_{1\times 1}& D_{1\times 2}&\cdots  & D_{1\times d}\\
D_{2\times 1}& D_{2\times 2} &\cdots  & D_{2\times d}\\
\cdots& \cdots&\cdots  & \cdots\\
D_{d\times1}& D_{d\times 2} &\cdots  & D_{d\times d}\\
\end{bmatrix}$

  $$ {\mathrm{where}} \quad D_{i\times j}= \begin{blockarray}{*{7}{c} l}
    \begin{block}{*{7}{>{$\footnotesize}c<{$}} l}
     &  &  & & & & \\
    \end{block}
    \begin{block}{[*{7}{c}]>{$\footnotesize}l<{$}}
     0& \cdot & \cdot & \cdot &\cdot & 0 &b_{0,j\times i} &\\
   0& \cdot &   \cdot & \cdot&\cdot & 0& 0 &  \\
   \cdot & \cdot & \cdot & \cdot &\cdot & \cdot & \cdot   &\\
       \cdot & \cdot & \cdot &\cdot & \cdot & \cdot & \cdot &\\
      0 & \cdot & \cdot & \cdot & \cdot &0 &  0 & \\
\end{block}
  \end{blockarray}
$$
\end{itemize}
Putting $\Psi=C$  in Theorem 8.1 of \cite{kk04} we obtain  the following exact sequence of \tm modules:
	$$0\lra \Phi^{\vee}\lra \Ext_{\tau}^1(\Phi,C)\lra \G_a^d\lra 0.$$
	However, further on we need the exact form of the \tm module  $\Phi^{\vee}$. This is the content of the next proposition. 
\begin{prop}\label{dual1}
Let $\Phi$ be a strictly pure \tm module given by (\ref{eqq}) then ${\Phi}^{\vee}$ has a structure of a \tm module. This structure is given by the matrix (\ref{mac}) and therefore
\begin{equation*}
\deg_\tau{\Phi}^{\vee}=2     \quad  {\mathrm{and}} \quad  \dim{\Phi}^{\vee}=d\cdot(n-1)
\end{equation*}
\end{prop}

\begin{proof}
We have ${\mathrm{Der}}(\Phi,C)={\mathrm{Mat}}_{1\times d}(L\{\tau\})$ (cf. \cite[last para. p.100]{kk04}).

Using formula (\ref{delta}) we determine the following inner biderivations:
\begin{align}\label{innerbi}
&{\delta}^{\big(c{\tau}^kE_{1\times i}A_n^{-1}\big)}(t)=c{\tau}^kE_{1\times i}A_n^{-1}\cdot\Phi_t- C_t\cdot c{\tau}^kE_{1\times i}A_n^{-1}=\\ \nonumber
=E_{1\times i}\cdot& \Bigg(c\Big({\theta}^{(k)}-\theta\Big)B_0^{(k)}{\tau}^k+\big(cB_1^{(k)}-c^{(1)}B_0^{(k+1)}\big){\tau}^{k+1}+\sum\limits_{j=2}^{n-1}B_j^{(k)}{\tau}^{k+j}\Bigg)\\\nonumber
+E_{1\times i}\cdot& c{\tau}^{n+j}
%{\delta}^{\big(c{\tau}^kE_{1\times i}A_n^{-1}\big)}=c{\tau}^kE_{1\times i}A_n^{-1}\Big(\big(\theta I_d+N_{\Phi}\big)\tau^0+\sum\limits_{i=1}^{\rk\Phi}A_i\tau^i\Big)-\\ \nonumber
%(\theta +\tau)c{\tau}^kE_{1\times i}A_n^{-1}=c\big({\theta}^{(k)}-\theta\big)E_{1\times i}B_0^{(k)}{\tau}^k+\Big(cE_{1\times i}B_1^{(k)}-\\\nonumber
%cE_{1\times i}B_0^{(k+1)}\Big){\tau}^{k+1}+\sum\limits_{j=2}^{n-1}B_j^{(k)}{\tau}^{k+j}+cE_{1\times i}{\tau}^{n+j},\quad k=0,1,\dots
\end{align}
 The inner biderivation
${\delta}^{\big(c{\tau}^kE_{1\times i}A_n^{-1}\big)}(t)$ at the  $i$-th coordinate  has degree $n+k$ and degree less than $n+k$ at all remaining coordinates. Hence,  proceeding analogously as in the proof of \cite[Lemma 8.1]{kk04} we obtain the following isomorphism: 
\begin{equation}\label{ext0}
\Ext_{0,\tau}^1(\Phi,C) \cong {\mathrm{Mat}}_{1\times d}\big(L\{\tau\}_{\langle 1, \deg_\tau\Phi )}\big).
\end{equation}
We choose the following basis:
\begin{equation}\label{basis}
\Big\{E_{1\times 1}{\tau}^k\Big\}_{k=1}^{n-1},\Big\{E_{1\times 2}{\tau}^k\Big\}_{k=1}^{n-1},\dots , \Big\{E_{1\times d}{\tau}^k\Big\}_{k=1}^{n-1},
\end{equation}
where
$\Big\{E_{1\times 1}{\tau}^k\Big\}_{k=1}^{n-1}$ 
is an abbreviation for $E_{1\times 1}{\tau}^1, E_{1\times 1}{\tau}^2, \dots, E_{1\times 1}{\tau}^{n-1}$. 
For the basis (\ref{basis}) we will find the \tm module structure on $\Ext_{0,\tau}^1(\Phi,C).$ To do this we have to compute the action of "t" on each of the basis elements and then 
perform the necessary reductions by inner derivations in order to   ensure that all representatives   belong  to the right-hand side of (\ref{ext0}). 
Notice that for a twisted polynomial $w(\tau)\in L\{\tau\}$ the expression
$w(\tau)_{|\tau=c}$
denotes the evaluation of this polynomial at $c.$ More generally for a vector of twisted polynomials we denote:
$\big[w_1(\tau),\dots, w_s(\tau)\big]_{\big|\tau=c}:= \big[w_1(\tau)_{|\tau=c},\dots, w_s(\tau)_{|\tau=c}\big].$

We have:
\begin{align}\label{tstr1}
t*E_{1\times 1}c{\tau}^k&=E_{1\times 1}\Big({\theta}c{\tau}^k+c^{(1)}{\tau}^{k+1}\Big)=\Big[\podwzorem{ 0,\dots, 0}{k-1}, \theta, \tau, 0,\dots, 0\Big]_{\big|\tau=c}\\\nonumber 
&\quad {\mathrm{for}} \quad k=1,\dots, n-2,
\end{align}
The right-hand sides of the second equalities of (\ref{tstr1}) represent coordinates of the vectors $t*E_{1\times 1}c{\tau}^k, \,\, k=1,\dots,n-2$ in the basis (\ref{basis}). 
Similarly, for $k=n-1$ we obtain: 
\begin{align}\label{tsrt2}
t*&E_{1\times 1}c{\tau}^{n-1}=E_{1\times 1}\Big({\theta}c{\tau}^{n-1}+\underbrace{c^{(1)}{\tau}^{n}}\Big)=\\ \nonumber
&=\Big[ \big(b_{0,1\times 1}^{(1)}c^{(2)}-b_{1,1\times 1}c^{(1)}\big)\tau-\sum\limits_{j=2}^{n-1}b_{j,1\times1}c^{(1)}{\tau}^j +{\theta}c{\tau}^{n-1},\\\nonumber
& \qquad\big(b_{0,1\times 2}^{(1)}c^{(2)}-b_{1,1\times 2}c^{(1)}\big)\tau
-\sum\limits_{j=2}^{n-1}b_{j,1\times2}c^{(1)}{\tau}^j,\dots,\\\nonumber
&\qquad \big(b_{0,1\times d}^{(1)}c^{(2)}-b_{1,1\times d}c^{(1)}\big)\tau
-\sum\limits_{j=2}^{n-1}b_{j,1\times d}c^{(1)}{\tau}^j\Big]=  \\\nonumber
&=\Big[b_{0,1\times 1}^{(1)}{\tau}^{2}-b_{1,1\times 1}{\tau},\big[-b_{j,1\times1}\tau\big]_{j=2}^{n-2},-b_{n-1,1\times 1}\tau+\theta, \\\nonumber
&\qquad b_{0,1\times 2}^{(1)}{\tau}^2-b_{1,1\times 2}\tau, \big[-b_{j,1\times2}\tau\big]_{j=2}^{n-1},\dots,\\\nonumber
&\qquad b_{0,1\times d}^{(1)}{\tau}^2-b_{1,1\times d}\tau, \big[-b_{j,1\times d}\tau\big]_{j=2}^{n-1}\Big]_{\big|\tau=c}.
\end{align} 
In formula  (\ref{tsrt2}) the under-braced element is reduced by the  inner biderivation ${\delta}^{\big(c^{(1)}{\tau}^0E_{1\times 1}A_n^{-1}\big)}.$
The right-hand side of the third equality of (\ref{tsrt2}) represents  coordinates of the vector $t*E_{1\times 1}c{\tau}^{n-1}$ in the basis (\ref{basis}).

Analogous computations yield the following formulas:
\begin{align}\label{tsrt3}
t*E_{1\times i}c{\tau}^k&=\Big[\podwzorem{\podwzorem{ 0,\dots, 0}{n-1},\dots, \podwzorem{ 0,\dots, 0}{n-1}}{i-1}, \podwzorem{ 0,\dots, 0}{k-1}, \theta, \tau, 0,\dots, 0\Big]_{\big|\tau=c}\\\nonumber 
&\quad {\mathrm{for}} \quad k=1,\dots, n-2,\\\label{tsrt4}
t*E_{1\times i}c{\tau}^{n-1}&=\Big[b_{0,i\times 1}^{(1)}{\tau}^{2}-b_{1,i\times 1}{\tau}, 
\big[-b_{j,i\times1}\tau\big]_{j=2}^{n-1},\dots, \\\nonumber
&\qquad\podwzorem{b_{0,i\times i}^{(1)}{\tau}^2-b_{1,i\times i}\tau, \big[-b_{j,i\times i}\tau\big]_{j=2}^{n-2}, \theta-b_{n-1, i\times i}\tau }{i-th \,\, block},\dots, \\\nonumber
&\qquad b_{0,i\times d}^{(1)}{\tau}^2-b_{1,i\times d}\tau, \big[-b_{j,i\times d}\tau\big]_{j=2}^{n-1}\Big]_{\big|\tau=c}.
\end{align}
Thus we obtain the following matrix for the dual \tm module ${\Phi}^{\vee}$:
\begin{equation}\label{mac}
{\Phi}_t^{\vee}={\theta}I_{d(n-1)}+A_1^{\vee}\tau +A_2^{\vee}{\tau}^2
\end{equation}
where $A_1^{\vee}=A_1^{\vee}(B_1,\dots,  B_{n-1})$ and  $A_2^{\vee}=A_2^{\vee}(B_0^{(1)})$. The proposition follows.
\end{proof}

The following proposition shows that the construction of a dual \tm module is functorial.
\begin{prop}\label{naturality}
Let $f:\Phi \rightarrow \Psi$ be a morphism of strictly pure \tm modules with no nilpotence of $\tau$-degree greater than  or equal to 2. Then there exists 
the following commutative diagram of \tm modules with the exact rows:
\begin{equation}\label{diagram 2.4}
\xymatrix{
0\ar[r] & {\Psi}^{\vee} \ar[d]^{f^{\vee}} \ar[r] & {\Ext}^1_{\tau}(\Psi,C) \ar[d]^{\bar f}\ar[r] & {\mathbb{G}}_a^{{\dim}{\Psi}}\ar[d]^{{\partial}f}\ar[r] & 0\\
0\ar[r] & {\Phi}^{\vee}  \ar[r]  & {\Ext}^1_{\tau}(\Phi,C) \ar[r] & {\mathbb{G}}_a^{{\dim}{\Phi}}\ar[r] & 0
}
\end{equation} 
where $f^{\vee}$ is the map of \tm modules induced by $f.$
\end{prop}
\begin{proof}
We begin with a definition of $\bar f$ on the level of biderivations. Let $\dim\Phi=d,$ $\dim\Psi=e,$ $\deg_\tau\Phi=n$ and $\deg_\tau\Psi=m.$
Let  $\Pi_t^{\Phi}$ and  $\Pi_t^{\Psi}$ be  matrices of the action of $t$ on
\begin{equation*}
 {\Ext}^1_{\tau}(\Phi,C) \cong  {\mathrm{Mat}}_{1\times d}(L\{{\tau}\})_{{\langle} 1, n)} \quad  \textnormal{and} \quad {\Ext}^1_{\tau}(\Psi,C) \cong  {\mathrm{Mat}}_{1\times e}(L\{{\tau}\})_{{\langle} 1, m)}, 
\end{equation*}
respectively.
Define
\begin{equation*}
{\bar f}\Big(\big[ {\delta}_1,\dots {\delta}_e\big]\Big):=\big[ {\delta}_1,\dots {\delta}_e\big]\cdot f
\end{equation*}
Since $f\in  {\mathrm{Mat}}_{e\times d}(L\{{\tau}\})$  we have $\big[ {\delta}_1,\dots {\delta}_e\big]\cdot f\in  {\mathrm{Mat}}_{1\times d}(L\{{\tau}\}).$
Assume that $\delta,\hat\delta\in  {\mathrm{Mat}}_{1\times e}(L\{{\tau\})}$ represent the same coset in ${\Ext}^1_{\tau}(\Phi,C)$ i.e. 
$\delta-\hat\delta = U\Psi-CU$ for some $U\in  {\mathrm{Mat}}_{1\times e}(L\{{\tau\})}.$ One readily verifies that
\begin{align*}
\bar{f}(\delta )-\bar{f}(\hat{\delta})=(\delta-\hat\delta )\cdot f=(U\Psi -CU)\cdot f=Uf{\Phi}-CUf\in \Der_{in}(\Phi,C).
\end{align*}
Thus $\bar f$ is well-defined if $f$ is a morphism of \tm modules. 
Moreover, since 
\begin{equation*}
\bar{f}(t*\delta)=C_t\delta\cdot f=C_t\cdot{\bar f}(\delta )=t*{\bar f}(\delta)
\end{equation*}
we see that $\bar f$ is a morphism of ${\mathbb F}_q[t]$-modules. Therefore ${\bar f}\cdot {\Pi}_t^{\Psi}={\Pi}_t^{\Phi}\cdot {\bar f}.$ It remains to be shown that in the bases 
for ${\Ext}^1_{\tau}(\Psi,C)$ and ${\Ext}^1_{\tau}(\Phi,C)$ the map  $\bar f $ is a map of \tm modules i.e. it is given by a matrix $F_0+F_1\tau \dots+F_k{\tau}^k\in  {\mathrm{Mat}}_{dn\times em}(L\{{\tau\})}$.
For this it is enough to find the image of $\bar f$ on the basis elements $E_{1\times l}\cdot c{\tau}^k$ for $l\in \{1,\dots e \}$ and $k\in \{0,\dots m-1 \}.$ If 
$\bar{f}(E_{1\times l}\cdot c{\tau}^k)$ has a coordinate with a degree of $\tau$ greater than $n-1$ then we reduce it  by means of the inner biderivation  (\ref{innerbi}). 
It is clear that we obtain the images of $\bar f$ on the basis elements as linear combinations of the basis elements of the corresponding basis elements of ${\Ext}^1_{\tau}(\Phi,C)$.
Therefore $\bar f$ is a morphism of \tm modules.
Define $f^{\vee}:={\bar f}|_{{\Psi}^{\vee}}.$ Since ${\partial}{\bar f}(\delta)\neq 0$ implies ${\partial}\delta\neq 0$ we obtain $f^{\vee}: {\Psi}^{\vee}\rightarrow {\Phi}^{\vee}.$ 
It is obvious that ${\partial}f :{\mathbb G}_{a}^{{\oplus}e} \rightarrow {\mathbb G}_{a}^{{\oplus}d} $ makes the diagram (\ref{diagram 2.4}) commute.
\end{proof}

Now we will determine $\Ext^1_{\tau}({\Phi}^{\vee},C)$ and ${\Phi}^{\vee\vee}=\Ext^1_{0,\tau}({\Phi}^{\vee},C)$. We need some more notation. 
Denote
\begin{equation*}
{\Phi}^{\vee}_t=\theta I_{d(n-1)}+A_1^{\vee}\tau + A_2^{\vee}{\tau}^2.
\end{equation*}
Then from Proposition \ref{dual1} one easily recovers $A_1^{\vee}$ and $A_2^{\vee}.$
Notice that the \tm module  $\Phi^{\vee}$ is not  strictly pure, therefore the structure of  t-action for $\Ext_{\tau}^1(\Phi^{\vee},C)$ cannot be obtained from  \cite[theorem 8.1]{kk04}. 
Moreover, ${\Phi}^{\vee}$ is not given by the composition series with the Drinfeld modules as  simple subquotients (cf. \cite{gkk}). 
We will show that despite these obstacles,  we can use the \tm reduction algorithm from  \cite{gkk},  but it will be much more complicated than in aforementioned  examples. 
The concept of the proof is related to that from
 \cite{pr}, where the corresponding result was proved for Drinfeld modules. However, in the multi-dimensional case a dual   \tm module is given by block matrices. This implies that  the corresponding biderivations, used for reductions, are in the block form. Therefore,  special care is needed in the reduction process.  We have to keep in mind that reductions in one block influence the others. 
This makes the situation much more complicated.
If we tried
to use   the inner biderivations of the form ${\delta}^{\big(c{\tau}^kE_{1\times l}\big)}$ in the \tm reduction algorithm, then we would face problems with  reductions in the columns $n-1, 2(n-1), \dots, d(n-1).$  The form of the matrix $A_2^{\vee}$   forces one to invert the matrix $B_0$ at each stage of reduction. As this matrix is  being changed in every reduction,  it is better to choose different inner biderivations for reductions. In order to describe these choices,  consider the following constructions. 

For the matrix $A_n=[a_{n,k\times l}]$  appearing in (\ref{eqq}), define the matrix $\widehat{A_n}$ as the following block matrix:
\begin{equation*}
\widehat{A_n}=\begin{bmatrix}
G_{1\times1} & G_{1\times 2} & \cdots & G_{1\times d}\\
G_{2\times1} & G_{2\times 2} & \cdots & G_{2\times d}\\
\cdot & \cdot & \cdots & \cdot\\
G_{d\times1} & G_{d\times 2} & \cdots & G_{d\times d}
\end{bmatrix}
\end{equation*}
where $[G_{k\times l}]$ are the following $(n-1)\times (n-1)$  matrices:
\begin{equation*}
 G_{k\times k}= \diag\big[\podwzorem{a_{n,k\times k}^{(1)},1,\cdots , 1}{n-1}\big],\quad k=1,\dots,d.
  \end{equation*} 
  and  
  \begin{equation*}
  G_{k,l}=a_{n,l\times k}^{(1)}\widehat{E}_{1\times(n-1)}\quad {\mathrm{for}}\quad k\neq l.
    \end{equation*} 
 $ \widehat{E}_{i\times j}$ denotes the $(n-1)\times(n-1)$ matrix with the only non-zero entry $1$ at the position $i\times j.$
One readily verifies that $\widehat{A_n}A_2^{\vee}=A_2^{\vee}(I_d)$ 
This  follows from the equality: $A_n^{(1)}B_0^{(1)}=I_d.$ 
We also compute
$\widehat{A_n}A_1^{\vee}=A_1^{\vee}(B_1A_n{^{(1)}},B_1,\dots, B_{n-1}).$
 Denote the elements of the matrix $B_1A_n^{(1)}$ as  $s_{1,i \times j}$  for $i,j =1,\dots d.$

We consider the following two types of inner biderivations:
\begin{align}\label{rdd1}
&{\delta}^{\big(c{\tau}^kE_{1\times (1+l(n-1))}\widehat{A}_n\big)}=E_{1\times (l+1)(n-1)}c{\tau}^{k+2} +\\ \nonumber +\Big[a_{n,1\times(l+1)}^{(k+1)}&c\big({\theta}^{(k)}-\theta\big)\tau^k-a_{n,1\times(l+1)}^{(k+2)}c^{(1)}{\tau}^{k+1},0,\dots,0, \podwzorem{-s_{1,1\times(l+1)}^{(k)}c{\tau}^{k+1}}{(n-1)-st \,\, column},\\\nonumber
a_{n,2\times(l+1)}^{(k+1)}&c\big({\theta}^{(k)}-\theta\big)\tau^k-a_{n,2\times(l+1)}^{(k+2)}c^{(1)}{\tau}^{k+1},0,\dots,0,       \podwzorem{-s_{1,2\times(l+1)}^{(k)}c{\tau}^{k+1}}{2(n-1)-st \,\, column},\dots 
 \\\nonumber 
 a_{n,d\times(l+1)}^{(k+1)}&c\big({\theta}^{(k)}-\theta\big)\tau^k-a_{n,d\times(l+1)}^{(k+2)}c^{(1)}{\tau}^{k+1},0,\dots,0,       \podwzorem{-s_{1,d\times(l+1)}^{(k)}c{\tau}^{k+1}}{d(n-1)-st \,\, column}\Big] 
 \\\nonumber 
 &\quad  {\mathrm{for}}\quad  l=0,\dots, (d-1)
\end{align}
and
\begin{align}\label{rdd2}
&{\delta}^{\big(c{\tau}^kE_{1\times (1+r+l(n-1))}\widehat{A}_n\big)}=\Big[0,\dots ,0, {c{\tau}^{k+1}}, \podwzorem{c\big({\theta}^{(k)}-\theta\big){\tau}^k-c^{(1)}{\tau}^{k+1}}{1+r+l(n-1)-st \,\, column} ,  0,\dots ,0
\Big]+\\\nonumber
  &+\Big[0,\dots ,0, \podwzorem{-b_{r+1,1\times(l+1)}^{(k)}}{(n-1)-st \,\, column},0, \dots,0 ,\podwzorem{-b_{r+1,2\times(l+1)}^{(k)} }{2(n-1)-th \,\, column},\dots , 
 \podwzorem{-b_{r+1,d\times(l+1)}^{(k)}}{d(n-1)-th \,\, column}  \Big]c{\tau}^{k+1}    \\\nonumber
  & {\mathrm{for}}\quad  l=0,\dots, (d-1) \quad  {\mathrm{and}}\quad r=1,\dots, (n-2)
\end{align}
\begin{rem}\label{rmrd1}
Notice that in (\ref{rdd2}) for $r=(n-2)$ the  term with the index $1+(n-2)+l(n-1)=(l+1)(n-1)$  equals $c\big({\theta}^{(k)}-\theta\big){\tau}^k-c^{(1)}{\tau}^{k+1}-
b_{n-1,(l+1)\times(l+1)}^{(k)}c{\tau}^{k+1}.$
\end{rem}
\begin{rem}\label{rmrd2}
 A matrix $U\in {\mathrm{Mat}}_{1\times d(n-1)}(L\{{\tau}\})$ can be divided into $d$ blocks  of size $1\times{(n-1)}.$ We enumerate these blocks in the following way:
 $U=\big[ U_{0},\dots , U_{d-1}\big]$. Then the inner biderivations (\ref{rdd1}) will be used for reductions of terms in the borders  of these blocks, i.e. in the columns
 $(n-1), 2(n-1),\dots , d(n-1),$  whereas the biderivations  (\ref{rdd2}) will be used for reductions inside the blocks, i.e.  for reductions of terms in  the remaining columns.
\end{rem}
Now we are ready to prove the following:
\begin{lem}\label{lemma}
We have:
\begin{itemize}
\item[(i)] $\Ext^1_{\tau}({\Phi}^{\vee},C)\cong \big[ \podwzorem{L\{{\tau}\}_{\langle 0,1)},\dots ,  L\{{\tau}\}_{\langle 0,1)}}{n-2}, L\{{\tau}\}_{\langle 0,2)}\big]^{\oplus d}
\\ \\\in  {\mathrm{Mat}}_{1\times d(n-1)}(L\{{\tau}\})$
\bigskip
\item[(ii)] $\Ext_{0,\tau}^1({\Phi}^{\vee},C)\cong \big[ \podwzorem{0,\dots ,  0 }{n-2}, L\{{\tau}\}_{\langle 1,2)}\big]^{\oplus d}
\in  {\mathrm{Mat}}_{1\times d(n-1)}(L\{{\tau}\})$
\end{itemize}
\end{lem}
\begin{proof}
Notice that (ii) follows directly from (i) and the fact that  ${\delta}^{\big(c{\tau}^0E_{1\times l}{\widehat A}_n\big)} \in \Der_{0,in}({\Phi}^{\vee},C)$ for $l=1,2,\dots , d(n-1).$
 So, we have to prove (i). 
 
 Let $\delta=\big[ {\delta}_l\big]_l=\big[ {\delta}_0,\dots {\delta}_d\big]\in \Der({\Phi}^{\vee},C).$ If the terms with the highest exponent of  $\tau$ are only on the borders  of the blocks, i.e. in the columns $(n-1), 2(n-1),\dots , d(n-1),$
 then we reduce them by the inner biderivations of type (\ref{rdd1}). In this way, we assure that the terms of the highest exponent of $\tau$ are inside the blocks.  
Now we will describe how to reduce the highest terms inside such a block. We start the reduction with the  highest leftmost term in this block. Let us say that this biderivation is in the 
$r$-th place of the $l$-th block, i.e. it has index $1\times l(n-1)+r$. We reduce this term by the inner biderivation ${\delta}^{\big(c_1{\tau}^kE_{1\times l(n-1)+r+1}{\widehat{A}}_n \big)}$ for some $c_1\in L$  
and  $k\geq 0.$ Because of the form of the derivation (\ref{rdd2}), we obtain the highest terms (of the same degree  we just reduced) at the border and also at the  next position to the  right of the reduced term.
Thus we "pushed" the highest terms in this block one position to the right. We continue reducing in this way until the highest terms are only at the borders of the blocks. 
Then we perform the necessary reductions by the inner biderivations of type (\ref{rdd1}). After this process our biderivation $\delta$ is reduced to a biderivation of degree less 
than that of $\delta.$ This describes the step of the downward induction with respect to the degree of the reduced biderivation. 
One easily sees that the induction ends if the reduced biderivation is on the right-hand side of (ii).
\end{proof}

\begin{thm}\label{thm2}
Let ${\Phi}$ be a strictly pure \tm module of $\tau$-degree $n>1$ and  dimension $d$ with no nilpotence and let ${\Phi}^{\vee}$ be its dual. Then ${\Phi}^{\vee\vee}\cong \Phi$ and there exists the following short exact sequence:
\begin{equation}\label{es}
0\rightarrow \Phi\rightarrow \Ext_{\tau}^1(\Phi^{\vee},C) \rightarrow \G_a^{(n-1)d}\rightarrow 0.
\end{equation}
\end{thm}

Before we prove Theorem \ref{thm2} it is worth 
giving an example that illustrates the idea of the proof.

\begin{ex}
    Let $\Phi=\theta I_2+\left[\begin{array}{cc}
      \alpha_1   & \alpha_2 \\
      \alpha_3   &  \alpha_4
    \end{array}\right]\tau+ 
    \left[\begin{array}{cc}
      \beta_1   & \beta_2 \\
      \beta_3   &  \beta_4
    \end{array}\right]\tau^2+
    \left[\begin{array}{cc}
      1   & 0 \\
      \gamma  &  1
    \end{array}\right]\tau^3$. 
    Applying notations from the begining of section 3 we obtain:
    \begin{align*}
        B_0 & =\left[\begin{array}{cc}
      1   & 0 \\
      \gamma  &  1
    \end{array}\right]^{-1}=\left[\begin{array}{cc}
      1   & 0 \\
      -\gamma  &  1
    \end{array}\right] \qquad
    B_2=
    \left[\begin{array}{cc}
      \alpha_1   & \alpha_2 \\
      \alpha_3-\gamma \alpha_1   &  \alpha_4-\gamma \alpha_2
    \end{array}\right]\\
    B_1&=
    \left[\begin{array}{cc}
      \beta_1   & \beta_2 \\
      \beta_3-\gamma \beta_1   &  \beta_4-\gamma \beta_2
    \end{array}\right]=:
    \left[\begin{array}{cc}
      \beta_1   & \beta_2 \\
      \widetilde{\beta}_3   &  \widetilde{\beta}_4
    \end{array}\right].
    \end{align*}
    Then according to the Proposition \ref{dual1} we get that $\Phi^\vee=\theta I_4+A_1^\vee\tau+A_2^\vee\tau^2,$ 
    $$A_1^\vee=\left[\begin{array}{cccc}
      0   & -\alpha_1 & 0 & - (\alpha_3-\gamma\alpha_1)  \\
      1   & -\beta_1 & 0 & - \widetilde{\beta}_3  \\
      0   & -\alpha_2 & 0 & - (\alpha_4-\gamma\alpha_2)  \\
      0   & -\beta_2 & 1 & - \widetilde{\beta}_4  \\
    \end{array}\right]\quad\textnormal{and}\quad 
    A_2^\vee=\left[\begin{array}{cccc}
      0   & 1 & 0 & - \gamma^{(1)}  \\
      0   & 0 & 0 & 0  \\
      0   & 0 & 0 & 1  \\
      0   & 0 & 0 & 0  \\
    \end{array}\right]
    $$
    We will show that $\Phi^{\vee\vee}\cong \Phi$. In order to do this, we first describe a \tm module structure on  $\Ext^1_{0,\tau}\big(\Phi^\vee, C\big)$.  
    By Lemma \ref{lemma} we obtain the following isomorphism:
    $$\Ext^1_{0,\tau}\big(\Phi^\vee, C\big)\cong \Big\{\big[0,a\tau, 0, b\tau\big]\mid a,b\in L\Big\}.$$
    Consider the following basis consisting of  the following two vectors: $[0,\tau,0,0]$ and $[0,0,0,\tau]$. In the reduction process we use the following inner biderivations described in Lemma  \ref{lemma}: 
    \begin{align*}
        \delta^{\big([c\tau^k,0,0,0]\cdot \widehat{A_3}\big)}, 
         \delta^{\big([0,c\tau^k,0,0]\cdot \widehat{A_3}\big)}, 
          \delta^{\big([0,0,c\tau^k,0]\cdot \widehat{A_3}\big)}, 
           \delta^{\big([0,0,0,c\tau^k]\cdot \widehat{A_3}\big)}, 
    \end{align*}
    where $\widehat{A_3}=\left[\begin{array}{cccc}
      1   & 0 & \gamma^{(1)} & 0  \\
      0   & 1 & 0 & 0  \\
      0   & 0 & 1 & 0  \\
      0   & 0 & 0 & 1 \\
    \end{array}\right].$ To ease notation let us denote
    $B_1\cdot \left[\begin{array}{cc}
      1   & 0 \\
      \gamma^{(1)}  &  1
    \end{array}\right]:=\left[\begin{array}{cc}
      s_1   & s_2 \\
      s_3  &  s_4
    \end{array}\right]$.\\
    %%%%%
    Now we compute the action of
     $t$ on the generators $[0,c\tau,0,0]$ and $[0,0,0,c\tau]$. Then
    \begin{align*}
        t&*\big[0,c\tau,0,0\big]= (\theta+\tau)\cdot \big[0,c\tau,0,0\big]= \Big[0,\theta c\tau+ c^{(1)}\tau^2,0,0\Big] \\
        &= \Big[c^{(2)}\tau,\big(\theta c + s_1c^{(1)}\big)\tau,c^{(2)}\gamma^{(2)}\tau, s_3c^{(1)}\tau\Big]\\
        &=
        \Big[0 ,\big(\theta c + s_1c^{(1)}+\beta_1c^{(2)}+c^{(3)}\big)\tau,c^{(2)}\gamma^{(2)}\tau, \big(s_3c^{(1)} + \widetilde{\beta}_3c^{(2)}\big)\tau\Big]\\
        &=\Big[0 ,\big(\theta c + s_1c^{(1)}+(\beta_1+\gamma^{(2)}\beta_2)c^{(2)}+c^{(3)}\big)\tau, 0, \big(s_3c^{(1)} + (\widetilde{\beta}_3+ \widetilde{\beta}_4 \gamma^{(2)})c^{(2)}+ c^{(3)}\gamma^{(3)}\big)\tau\Big]\\
        &=\Big[\theta  + s_1\tau+(\beta_1+\gamma^{(2)}\beta_2) \tau^2+\tau^3, s_3\tau + (\widetilde{\beta}_3+ \widetilde{\beta}_4 \gamma^{(2)})\tau^2+ \gamma^{(3)}\tau^3\Big]_{\big|\tau=c},
    \end{align*}
    where the equalities $3,4$ and $5$ result  from the reductions by means of the biderivations  
     $\delta^{\big([c^{(1)}\tau^0,0,0,0]\cdot \widehat{A_3}\big)}, $
     $\delta^{\big([0,c^{(2)}\tau^0,0,0]\cdot \widehat{A_3}\big)} $
     and $\delta^{\big([0,0,0,c^{(2)}\gamma^{(2)}\tau^0]\cdot \widehat{A_3}\big)}$ respectively.\\ 
     Similarly, 
     \begin{align*}
    t&*\big[0,0,0,c\tau\big]= \Big[ s_2\tau +\beta_2\tau^2, \theta + s_4\tau+\widetilde{\beta}_4 \tau^4 +\tau^3 \Big]_{\big|\tau=c}.  
     \end{align*}
     Thus,
     $$\Phi^{\vee\vee}=\theta I_2+\left[\begin{array}{cc}
      s_1   & s_2 \\
      s_3   &  s_4
    \end{array}\right]\tau+ 
    \left[\begin{array}{cc}
      \beta_1+\beta_2\gamma^{(2)}   & \beta_2 \\
      \widetilde{\beta}_3+\widetilde{\beta}_4\gamma^{(2)}   & \widetilde{\beta}_4
    \end{array}\right]\tau^2+
    \left[\begin{array}{cc}
      1   & 0 \\
      \gamma^{(3)}  &  1
    \end{array}\right]\tau^3.$$
    Straightforward calculation yields the following equality
    $$\Phi^{\vee\vee}=\left[\begin{array}{cc}
      1   & 0 \\
      \gamma  &  1
    \end{array}\right]^{-1}\cdot \Phi\cdot \left[\begin{array}{cc}
      1   & 0 \\
      \gamma  &  1
    \end{array}\right],$$
  and therefore   an isomorphism $\Phi^{\vee\vee}\cong \Phi.$
\end{ex}

\begin{proof}[Proof of  Theorem \ref{thm2}]
First we prove the existence of the exact sequence (\ref{es}). From lemma \ref{lemma} it follows that the basis: 
\begin{equation*}
\Big\{ E_{1\times l(n-1)+2}{{\widehat A}_n}, E_{1\times l(n-1)+3}{{\widehat A}_n}, \dots ,  E_{1\times (l+1)(n-1), }{{\widehat A}_n},  E_{1\times l(n-1)+1}{{\widehat A}_n}  \Big\}_{l=0}^{d-1}
\end{equation*}
satisfies step 2 of the reduction algorithm of \cite{gkk}. Therefore from \cite[Proposition 3.1]{gkk} the $\F_q[t]$-module $\Ext_{\tau}^1(\Phi^{\vee},C)$ is a \tm module of dimension $n\cdot d.$ 
Moreover, every inner biderivation ${\delta}^{(c{\tau}^0E_{1\times l}{\widehat A}_n)} \in \Der_{0,in}({\Phi}^{\vee},C)$ and therefore the \tm module structure on $\Ext_{0,\tau}^1(\Phi^{\vee},C)$ is induced from the \tm module structure of   $\Ext_{\tau}^1(\Phi^{\vee},C)$ 
by deleting $(n-1)d$ rows and  $(n-1)d$ columns that correspond to $E_{1\times i}{\tau}^0,\,\, i=1,\dots d(n-1).$ This shows (\ref{es}).
It remains to be shown that ${\Phi}^{\vee\vee}\cong{\Phi}.$ In order to do this we have to determine the exact \tm module structure on  $\Ext_{0,\tau}^1(\Phi^{\vee},C).$

Since $\Ext_{0,\tau}^1(\Phi^{\vee},C) \cong \big[ \podwzorem{0,\dots , 0}{n-2},L\{\tau\}_{\langle 1,2)}      \big]^{{\oplus}d}$
we choose the following basis:  
\begin{equation}\label{bbbs} 
E_{1\times {l+1}(n-1)}\tau \quad \textnormal{for}\quad l=0,\dots d-1. 
\end{equation}
Then 
\begin{equation*}
t*E_{1\times (l+1)(n-1)}c\tau=(\theta +\tau)cE_{1\times (l+1)(n-1)}=E_{1\times (l+1)(n-1)}\big({\theta}c\tau+c^{(1)}{\tau}^2\big)
\end{equation*}
We reduce the term $E_{1\times (l+1)(n-1)}c^{(1)}{\tau}^2$ by the  biderivation  ${\delta}^{\big(c^{(1)}{\tau}^0E_{1\times l(n-1)}\widehat{A}_n\big)}$  of the type (\ref{rdd1}):
and after the reduction we obtain:
\begin{align*}
t*E_{1\times (l+1)(n-1)}c\tau=\Big[&a_{n,1\times(l+1)}^{(2)}c^{(2)}c{\tau},0,\dots ,0, s_{1,1\times(l+1)}c^{(1)}{\tau},\dots,\\ &\podwzorem{a_{n,(l+1)\times(l+1)}^{(2)}c^{(2)}\tau,0,\dots,0, \theta c\tau+ s_{1,(l+1)\times(l+1)}c^{(1)}{\tau}}{l-th\quad block},\dots \\ 
&a_{n,d\times(l+1)}^{(2)}c^{(2)}\tau,0,\dots,0, s_{1,d\times(l+1)}c^{(1)}{\tau} \Big] 
\end{align*}
Now we describe how to reduce an element inside the $m-$th block.
 We use, as in the proof of lemma \ref{lemma}, the following sequence of inner biderivations of type (\ref{rdd2}):
${\delta}^{\big(c_s \tau^0E_{1\times m(n-1)+s}\widehat{A}_n\big) }$ for $ s=2,3,\dots, n-1,$   
where $c_s= a_{n,(m+1)\times (l+1)}^{(s)} c^{(s)}$. 
After applying this sequence of inner biderivations for each block from $m=0,1,\dots d-1$ we obtain the following reduced form of the product:
\begin{align*}
&t*E_{1\times (l+1)(n-1)}c\tau=E_{1\times (n-1)}\Big(s_{1, 1\times (l+1)}c^{(1)}+\sum\limits_{m=0}^{d-1}b_{2,1\times(m+1)}a_{n,(m+1)\times(l+1)}^{(2)}c^{(2)}\\\nonumber
& +\sum\limits_{m=0}^{d-1}b_{3,1\times(m+1)}a_{n,(m+1)\times(l+1)}^{(3)}c^{(3)}+ \dots   +\sum\limits_{m=0}^{d-1}b_{n-1,1\times(m+1)}a_{n,(m+1)\times(l+1)}^{(n-1)}c^{(n-1)}+ \\\nonumber
& +a_{n,1\times(l+1)}^{(n)}c^{(n)}\Big)\tau+E_{1\times 2(n-1)}\Big(s_{1, 2\times (l+1)}c^{(1)}+\sum\limits_{m=0}^{d-1}b_{2,2\times(m+1)}a_{n,(m+1)\times(l+1)}^{(2)}c^{(2)}\\
& +\sum\limits_{m=0}^{d-1}b_{3,2\times(m+1)}a_{n,(m+1)\times(l+1)}^{(3)}c^{(3)}+ \dots   +\sum\limits_{m=0}^{d-1}b_{n-1,2\times(m+1)}a_{n,(m+1)\times(l+1)}^{(n-1)}c^{(n-1)}+\\
\end{align*}
\begin{align*}
& +a_{n,2\times(l+1)}^{(n)}c^{(n)}\Big)\tau+\cdots+\\
&\qquad \qquad \qquad\qquad +E_{1\times d(n-1)}\Big(s_{1, d\times (l+1)}c^{(1)}+\sum\limits_{m=0}^{d-1}b_{2,d\times(m+1)}a_{n,(m+1)\times(l+1)}^{(2)}c^{(2)}\\\nonumber
& +\sum\limits_{m=0}^{d-1}b_{3,d\times(m+1)}a_{n,(m+1)\times(l+1)}^{(3)}c^{(3)}+ \dots   +\sum\limits_{m=0}^{d-1}b_{n-1,d\times(m+1)}a_{n,(m+1)\times(l+1)}^{(n-1)}c^{(n-1)}+\\
&+a_{n,d\times(l+1)}^{(n)}c^{(n)}\Big)\tau
 \end{align*}
Thus in the basis (\ref{bbbs})  we obtain the following expansion:
\begin{align}\label{final2}
t*E_{1\times (l+1)(n-1)}c\tau=\Big[s_{1, 1\times (l+1)}\tau+\sum\limits_{m=0}^{d-1}b_{2,1\times(m+1)}a_{n,(m+1)\times(l+1)}^{(2)}{\tau}^{2}\\\nonumber
 +\sum\limits_{m=0}^{d-1}b_{3,1\times(m+1)}a_{n,(m+1)\times(l+1)}^{(3)}{\tau}^{3}+ \dots   +\sum\limits_{m=0}^{d-1}b_{n-1,1\times(m+1)}a_{n,(m+1)\times(l+1)}^{(n-1)}{\tau}^{n-1}
 \\\nonumber
+a_{n,1\times(l+1)}^{(n)}{\tau}^{n} ,s_{1, 2\times (l+1)}\tau+\sum\limits_{m=0}^{d-1}b_{2,2\times(m+1)}a_{n,(m+1)\times(l+1)}^{(2)}{\tau}^{2}\\\nonumber
 +\sum\limits_{m=0}^{d-1}b_{3,2\times(m+1)}a_{n,(m+1)\times(l+1)}^{(3)}{\tau}^{3}+ \dots   +\sum\limits_{m=0}^{d-1}b_{n-1,2\times(m+1)}a_{n,(m+1)\times(l+1)}^{(n-1)}{\tau}^{n-1}\\
\nonumber
+a_{n,2\times(l+1)}^{(n)}{\tau}^{n},\dots
 ,s_{1, d\times (l+1)}\tau+\sum\limits_{m=0}^{d-1}b_{2,d\times(m+1)}a_{n,(m+1)\times(l+1)}^{(2)}{\tau}^{2}\\\nonumber
 +\sum\limits_{m=0}^{d-1}b_{3,d\times(m+1)}a_{n,(m+1)\times(l+1)}^{(3)}{\tau}^{3}+ \dots   +\sum\limits_{m=0}^{d-1}b_{n-1,d\times(m+1)}a_{n,(m+1)\times(l+1)}^{(n-1)}{\tau}^{n-1}\\ \nonumber
 +a_{n,d\times(l+1)}^{(n)}{\tau}^{n}\Big]_{\big|\tau=c}
 \end{align}
But formula (\ref{final2}) written in  matrix form  yields the following \tm module structure on 
$\Ext_{0,\tau}^1(\Phi^{\vee},C)$:
\begin{equation}\label{final3}
{\Phi}^{\vee\vee}_t=I_d\theta +B_1A_n^{(1)}\tau+B_2A_n^{(2)}{\tau}^2+\dots +B_{n-1}A_n^{(n-1)}+A_n^{(n)}{\tau}^n=A_n^{-1}{\Phi}A_n.
\end{equation}
Therefore we obtained the isomorphism ${\Phi}^{\vee\vee}\cong {\Phi}.$

\end{proof}

Now we are ready to prove Theorem \ref{main}. 
\begin{proof}
Point a) is a specialization  of  \cite[Proposition 8.2]{kk04}. Point b) follows from Theorem \ref{thm2}. Finally point c) follows from Proposition \ref{naturality}.
\end{proof}

\section{Concluding remarks}
In section 3 we proved the Weil-Barsotti formula for strictly pure \tm modules with no nilpotence.  There are many theorems which one can  extend from the category
of Drinfeld modules to the subcategory (of the category of \tm modules) of strictly  pure \tm modules (with non-zero nilpotent part).  Good examples of this analogy are various theorems concerning heights (see \cite{de}, \cite{po}).
However, the following example shows that, concerning duality, we have to assume that the \tm module ${\Phi}_t$ is not only strictly pure but also that it has no nilpotence.
\begin{example}\label{exs}
Consider the following \tm module:
\begin{equation}\label{exs1}
{\Phi}_t=\theta I_3+\begin{bmatrix}
0&0&0\\
0&0&0\\
a&0&0
\end{bmatrix}
+I_3{\tau}^3 \quad {\mathrm{for}} \quad a\neq 0.
\end{equation}
Then applying the procedure used in the proof of Proposition  \ref{dual1} we obtain:
\begin{equation}\label{exs2}
{\Phi}_t^{\vee}=\begin{bmatrix}
\theta&{\tau}^2&0&0&0&0&-a^{(1)}\tau\\
\tau&\theta&0&0&0&0&0\\
0&0&\theta&{\tau}^2&0&0&0\\
0&0&\tau&\theta&0&0&0\\
0&0&0&0&\theta&0&(\theta-{\theta}^{(1)})\tau\\
0&0&0&0&\tau&\theta&{\tau}^2\\
0&0&0&0&0&\tau&\theta\\
\end{bmatrix}
\end{equation}
Notice that ${\Phi}_t^{\vee}$  has the zero nilpotent part but is not strictly pure.
Notice also that in order to apply the \tm reduction algorithm we have to assume that $L$ is a perfect field.
Then again using  the procedure described in the proof of theorem \ref{thm2} we obtain:
\begin{equation}\label{exs3}
\Ext^1_0({\Phi}_t^{\vee},C)_t=
\begin{bmatrix}
\theta+{\tau}^3&0&0\\
0&\theta+{\tau}^3&0\\
0&0&\theta+{\tau}^3-1
\end{bmatrix}
\end{equation}
We see that (\ref{exs3} ) is not a \tm module.
\end{example}
\begin{rem}\label{last}
Example \ref{exs} shows that for a \tm module which is not strictly pure,  its dual might not exist in the category of \tm modules.
However,  applying duality to our isomorphism ${\Phi}^{\vee\vee}\cong {\Phi}$ for ${\Psi}={\Phi}^{\vee}$ we get ${\Psi}^{\vee\vee}\cong {\Psi}.$
But $\Psi$ is not necessarily a strictly pure \tm module. This was the reason we had to use the very  special order of reductions in the proof of Theorem \ref{thm2}.
So we see that the assumptions on ${\Phi}$ in Theorem {\ref{main}} are sufficient  but not necessary.
\end{rem}
\section{Acknowledgment}
The authors would like to thank the  referee for numerous corrections,
valuable questions, comments and suggestions.
\appendix
\section{}
In this appendix we relate our construction to the duality considered in \cite{ta}.
For the sake of simplicity we will keep the same notation as in \cite{ta}. We advise the reader to consult \cite{ta} for the basic notions of a $\phi$-module, finite \tm module etc.  

 Recall that Y. Taguchi in \cite{ta} considers duality in the context of {\bf finite} \tm modules. For a Drinfeld module $\phi$ this amounts to the construction of the Weil pairing $_{a}\phi\times _{a}{\phi}^*\rightarrow _{a}C,$ where $a\in{\mathbb F}_q[t], $ $_{a}\phi$ is a finite \tm module of $a$-torsion points of a Drinfeld module, $_{a}{\phi}^*$ is the dual to the finite \tm module  $_{a}\phi$ and $_aC$ is the finite \tm module of $a$-torsion points of the Carlitz module. For a Drinfeld module Y. Taguchi writes the form of what he calls the dual Drinfeld module $\phi^\vee$ and justifies this name by showing that $_a\phi^\vee\cong\, _a{\phi}^*$ for all $a\in\F_q[t]$. In \cite{pr} it is shown that $\phi^\vee$ is in fact equal to $\Ext^1_{\tau,0}(\phi,C).$ The author of \cite{ta} remarks that one can construct similar Weil pairing 
 $_a\Phi\times_R \,_a\Phi^*\rightarrow _aC$ for $\Phi$ a strictly pure \tm module with no nilpotence. Since we have computed the form of $\Phi^\vee=\Ext^1_{\tau,0}(\Phi,C)$ we will prove that $_a\Phi^\vee\cong _a{\Phi}^*.$
 We have the following :
 \begin{lem}\label{aq}
 Assume $\Phi_t=\theta I+A_1X^q+\dots A_nX^{q^n},$ $A_i\in M_{d\times d}(R),$ $A_n\in \GL_d(R),$ and $X=(X_1,\dots , X_d)^T$ is a strictly pure \tm module with respect to the trivialization $E\cong \G_a^d=\Spec R[X_1,\dots, X_d]$  then for $a\in {\mathbb F}_q[t],$ $_a\Phi:=\mathrm{Ker}(\Phi_a)$ is a finite \tm module. This finite \tm module can be endowed with the standard $v$-module structure.
 \end{lem}
 \begin{proof}
  For an algebraic group $G$ over $S$ let ${\cal E}_G:={\underline\Hom}_{{\mathbb F}_q,S}(G,\G_a),$ where $S=\Spec R$ and $\G_a$ is the additive algebraic group over $S,$ be the 
  Zariski sheaf on $S$ of the ${\mathbb F}_q$- linear homomorphisms.
   Let  ${\cal E}_G^{(q)}:={\cal E}_G\otimes_{{\cal O}_S}{\cal O}_S$ be the base change of ${\cal E}$ by the $q$-th power map ${\cal O}_S\rightarrow {\cal O}_S.$ 
  Now, let $G:= {\mathrm{Ker}}({\Phi}_a).$ Then ${\cal E}_G$ is a free  $R$-module of rank 
  $n\cdot d\cdot \deg (a) $ with a basis $\{X_i^{q^j}: 1\leq i\leq d,\,\, 0\leq j\leq n\cdot\deg(a)-1\}.$ 
   We also see that $\phi: {\cal E}_G^{(q)}\rightarrow {\cal E}_G$ is given by the following formula:
  $$\phi\big(X_i^{q^j}\otimes 1\big)= X_i^{q^{j+1}}.$$
  Notice that $X_i^{q^{j+1}}$ for $j+1\geq n\cdot\deg(a)$ has to be calculated by means of the relation ${\Phi}_a(X)=0, \quad X=(X_1,\dots ,X_d )^T.$ This is possible since the matrix $A_n$ is invertible.

  As far as the $v-$structure is concerned, we need to find $v:{\cal E}_{\mathrm{Ker}{\Phi_a}}\rightarrow {\cal E}^{(q)}_{\mathrm{Ker}{\Phi_a}}$ such that $\Phi_t=\theta I+\phi\circ v$ (cf. \cite[Def.(3.2)]{ta}).
  We will define $v$ on  the basis $\big\{X_i^{q^j}: 1\leq i\leq d,\,\, 0\leq j\leq n\cdot\deg(a)-1\big\}.$
  Assume that for $i=1,\dots , d $ we have
  \begin{align}
  \Phi_t(X_i)=\theta IX_i+\sum_{k=1}^nA_kX_i^{(q^k)}=
  \theta X_i+\sum_{k=1}^n \sum_{l=1}^da_{k,l,i}X_i^{(q^k)}
  \end{align}
  Then we put (cf. \cite[Example (3.4)]{ta})
 \begin{equation}\label{vstr}
  v\big(X_i^{q^j}\big)=X_i^{q^{j-1}}\otimes \big({\theta}^{q^j}-\theta \big)+\sum_{k=1}^n \sum_{l=1}^dX_i^{q^{j+k-1}}\otimes a_{k,l,i}^ {q^j}
  \end{equation}
 \end{proof}

 For the sake of completeness we will prove the following, analogous to \cite[Theorem 5.1]{ta} theorem, 
 where the determined by us form  \eqref{mac}
 of $\Phi^\vee$ is necessary.
\begin{thm}Let $\Phi$ be a strictly pure \tm module with no nilpotence of the form (\ref{eqq}).
\begin{itemize}
    \item[$(i)$] If R is a perfect field, then   
$\Phi^\vee$ is an abelian \tm module of
$t$-rank $r(\Phi^\vee)= d\cdot n,$ $\tau$-rank $\rho(\Phi^\vee)= d\cdot (n-1),$ and weight $w(\Phi^\vee)= (n-1)/n$ in the
sense of \cite{a}.
    \item[$(ii)$]  For a non-zero $a\in A,$ the kernel $_a\Phi^\vee$ of the action of a on $\Phi^\vee$ is a finite \tm module over $R$ of rank $q^{r·deg(a)}.$
    \item[$(iii)$] For a non-zero $a \in A,$ there exists an $A$-bilinear pairing defined over $R$:
$$_a\Pi_\Phi :\,  _a\Phi\times_R \,_a\Phi^\vee \rightarrow \,\, _aC.$$
\item[$(iv)$] If we furnish $_a\Phi$ with the  $v$-module structure (cf. \eqref{vstr}) then
we have $_a\Phi^\vee\cong _a\Phi^*$
 and the pairing $_a\Pi_\Phi$ of (iii) coincides with the pairing $\Pi_\Phi$ of \cite[Theorem 4.3]{ta}.
\end{itemize}
\end{thm}
\begin{proof}
The idea of the proof is the same as that of \cite[Theorem 5.1]{ta} although the details are naturally different because our dual ${\Phi}^\vee$ is now more complicated. 

(i)
We need to find a suitable $R[\Phi_t]$ -basis of $\End_{{\mathbb F}_q,R}\big({\mathbb G}_a^{d\cdot(n-1)},{\mathbb G}_a\big).$ We will use the following basis: $$\big\{A_n^{-1}Y_{i,(n-1)}^q, Y_{i,1},\dots, Y_{i,(n-1)},\quad i=1,\dots , d\big\}. $$ This directly follows from the formulas 
(\ref{mac}) and notations adapted at the beginning 
of section \ref{Main results}. Note that $B_0=A_n^{-1}.$ Thus $r=d\cdot n,$ $\rho =d\cdot (n-1)$ and $w=\frac{n-1}{n}.$

(ii) Let $G=_a\Phi^\vee.$ The affine ring ${\cal O}_G$  has the following form:
$$R[Y_{1,1},\cdots,Y_{1,(n-1)},\cdots , Y_{d,1},\cdots,Y_{d,(n-1)}]/\Phi^\vee_a ({\mathbb Y})$$
 We have  to show that ${\cal O}_G$ is free over $R$ of rank $q^{r\cdot\deg (a)}$ and also that ${\cal E}_G$ is free over $R$ of rank $r\cdot \deg (a).$ So, let $a\in A={\mathbb F}_q[t]$ be a monic polynomial of degree $k$: $a=t^k+{\sum}_{l=0}^{k-1}g_l t^l.$

Define $Y_{s,i,j}\in {\cal O}_G$ for $1\leq s\leq d,\,\,0\leq i \leq k-1,\,\, 1\leq j\leq n-1$ by the following recursion: 
Similarly as in \cite{ta}, slightly abusing notation, we set:
\begin{equation}\label{rec11}
Y_{s,k-1,j}=Y_{s,j} \qquad 1\leq j\leq n-1,
\end{equation}
and  $${\mathbb Y}_{s,k-1}=(\widehat{Y}_{s,k-1,1},\dots, \widehat{Y}_{s,k-1,\rho})^T$$ where
$$\widehat{Y}_{s,k-1,l}=\begin{cases}
Y_{s,k-1,l}  \quad {\mathrm{for}} \quad  l=(s-1)\cdot (n-1)+j, \quad 1\leq j \leq n-1 \\
0 \quad {\mathrm{otherwise}}
    \end{cases}.$$
Define
\begin{equation}\label{af1}
{\mathbb Y}_{s,i-1} =\Phi_t^\vee ({\mathbb Y}_{s,i})+g_i {\mathbb Y}_{s,k-1}, \quad 1\leq i\leq k-1,  \quad 1\leq s\leq d
\end{equation}
where ${\mathbb Y}_{s,i}:=(\widehat{Y}_{s,i,1},\dots, \widehat{Y}_{s,i,\rho})^T. $
Then we see that
$${\mathbb Y}_{s,0}={\Phi}^\vee _{(t^{k-1}+g_{k-1}t^{k-2}+\cdots+g_1)}({\mathbb Y}_{s,k-1})$$ and ${\Phi}^\vee _t ({\mathbb Y}_{s,0})={\Phi}^\vee _a({\mathbb Y}_{s,k-1})-g_0{\mathbb Y}_{s,k-1}.$

Thus $G=_a\Phi^\vee$ is equivalent to $\Phi^\vee_a({\mathbb Y}_{k-1})=0.$ This yields the following equality: 
\begin{equation}\label{af2}
\Phi_t({\mathbb Y}_{s,0})=-g_0{\mathbb Y}_{s,k-1}
\end{equation}
Therefore one can consider $\cal O_G$ as the quotient of the following ring:
\begin{equation}\label{fa3}
R\big[Y_{s,i,j}, \,\, 1\leq s \leq d,\,\, 0\leq i\leq k-1, \,\, 1\leq j\leq n-1]
\end{equation}
by the relations (\ref{af1}) and (\ref{af2}).
Now, we modify an elegant argument of Taguchi. We put  $Y^{\prime}_{s,i,j}=Y_{s,i,j}^q$ for $j<n-1$ and $Y^{\prime}_{s,i,(n-1)}=Y_{s,i,(n-1)}$ and we can regard ${\cal O}_G$ as the quotient ${\cal O}^{\prime}$ of the ring
\begin{equation}\label{fa4}
 R\big[Y^{\prime}_{s,i,j}, \,\, 1\leq s \leq d,\,\, 0\leq i\leq k-1, \,\, 1\leq j\leq n-1]   
\end{equation}
by the same relations.
The key point is that in this ring the relations 
(\ref{af1}) and (\ref{af2}) can be written in the following form:
\begin{equation}\label{fa5}
 Y^{\prime q^2}_{s,i,j} +{\mathrm{lower\,\,  terms}}=0, \quad 1\leq s\leq d,\, \, 0\leq i\leq k-1,\,\, 1\leq j\leq n-1
\end{equation}
In our case this follows from the formula (\ref{mac}) for $\Phi_t^\vee$ and the form of block matrices $C_{i,j}$ and $D_{i,j}$ for $i,j=1,\dots d.$ One can solve the corresponding system of linear equations coming from (\ref{af1}) and (\ref{af2}) for the highest terms.
By \cite[Lemma 1.9.1]{AT} ${\cal O}^{\prime}$ is free over $R$ of rank $q^{2k\cdot d\cdot (n-1)}.$
By (\ref{fa5}) we see that the following set:
$${\prod}_{s,i,j}(Y^{\prime }_{s,i,j})^{l_{s,i,j}}; \quad 0\leq l_{s,i,j}\leq q^2-1.$$
${\cal O}_G$ is then a free submodule of ${\cal O}^{\prime}$ with the following basis:
$${\prod}_{s,i,j}(Y^{\prime }_{s,i,j})^{l_{s,i,j}}; \quad q|{l_{s,i,j}} \,\, {\mathrm{if}}\,\, 1\leq j\leq n-2.$$
Thus the rank of ${\cal O}_G$ equals $q^{kd(n-2)}\cdot q^{2kd}=q^{kdn}=q^{kr}=q^{\deg(a)\cdot r}.$ On the other hand ${\cal E}_G$ is free over $R$ with the following basis:
$$\{ Y_{s,i,j}\,\, 1\leq s\leq d, \,\,0\leq i\leq k-1,\,\, 0\leq j\leq n-1\} .$$ Therefore $\mathrm{rank}({\cal O}_G)=q^{\mathrm{rank}({\cal E}_G)}.$

(iii) Passing to affine rings one has to construct a map:
$${\cal O}_{_aC}\rightarrow {\cal O}_{_a\Phi}\otimes {\cal O}_{_a\Phi^\vee}$$
\begin{align*}\label{pi}
\pi: R[Z]/{\gamma}_a(Z)\rightarrow R[X_1,\dots ,X_d]/\Phi_a({\mathbb X}) \otimes \\R[Y_{s,i}, \,\, 1\leq s\leq d,\,\, 1\leq i\leq (n-1)]/\Phi^\vee(\mathbb Y),
\end{align*}
where $\gamma_t(Z)=\theta Z+Z^q$ is the Carlitz module $C,$ such that the following diagram:
\begin{equation}\label{diagram pi}
\xymatrix{
{\cal O}_{_aC}\ar[r]^{\pi\qquad}\ar[d]_{\gamma_t} & {\cal O}_{_a\Phi}\otimes {\cal O}_{_a\Phi^\vee}\ar[d]^{\Phi_t \otimes 1,\, 1\otimes\Phi_t^\vee}\\
{\cal O}_{_aC}\ar[r]^{\pi\qquad} & {\cal O}_{_a\Phi}\otimes {\cal O}_{_a\Phi^\vee}
}
\end{equation} commutes.
Let  $a=t^k+\sum_{i=0}^{k-1}g_it^i$ and $Y_{s,i,j}\in {\cal O}_{_a\Phi^\vee}$ be as in the proof of (ii).
Let $Y_{s,i,0}:=A_n^{-1}Y^q_{s,i,k-1},$ $1\leq s\leq d,$ $0\leq i\leq k-1.$ Denote $X_{s,i,j}:=\Phi_{t^i}(X_s)^{q^j}.$ Thus
\begin{equation}\label{fa14}
X_{d,i+1,0}=\Phi_t(X_{s,i,0})=\theta X_{s,i,0}+\sum_{j=1}^nA_jX_{s,i,j}
\end{equation}
and
\begin{equation}\label{fa15}
 0=\Phi_a(X_s) =X_{s,k,0}+\sum_{i=0}^{k-1}g_iX_{s,i,0}.  
\end{equation}
From (\ref{fa14}) and(\ref{fa15}) we get:
\begin{equation}\label{fa16}
0=\theta X_{s,k-1,0}+\sum_{j=1}^nA_jX_{s,k-1,j} +\sum_{i=0}^{k-1}g_iX_{s,i,0}
    \end{equation}
Define
\begin{equation}\label{fa15a}
\pi : Z\rightarrow \sum_{s=1}^d\sum_{i=0}^{k-1}\sum_{j=0}^{n-1}X_{s,i,j}\otimes Y_{s,i,j}
\end{equation}
The map $\pi$ is compatible with co-multiplications and the multiplication by ${\mathbb F}_q.$ One has to show commutativity of the diagram (\ref{diagram pi}).
By (\ref{fa15}) raised to the $q$-th power we obtain the following:
\begin{equation}\label{fa16}
 (\Phi_t\otimes 1)\pi(Z)=\sum_{s=1}^d\sum_{i=0}^{k-1}\sum_{j=0}^{n-1}X_{s,i+1,j}\otimes Y_{s,i,j}   \end{equation}
 \begin{align*}
 =-\sum_{s=1}^d \sum_{j=0}^{n-1}X_{s,0,j}\otimes g_0Y_{s,k-1,j}+\sum_{s=1}^d\sum_{i=0}^{k-1}\sum_{j=0}^{n-1}X_{s,i,j}\otimes (Y_{s,i-1,j}-g_iY_{s,k-1,j})
 \end{align*}
 By (\ref{af1}) and (\ref{af2}) we obtain  the following equality:
 $$(\Phi_t^\vee\otimes 1)\pi=(1\otimes\Phi_t^\vee )\pi$$
 Calculation of $\pi\gamma_t(Z)$ follows  the Taguchi's proof as well - one has to add another summation and replace scalars $a_j, 1\leq j\leq n$ by the corresponding matrices $A_j.$ We leave this calculation to  the reader.

 (iv) ${\cal E}_{_a\Phi}$ and ${\cal E}_{_a\Phi^\vee}$ may by viewed as dual to each other once making $X_{s,i,j}$ and $Y_{s,i,j}$ the dual bases. To see that $_a\Phi^\vee\cong _a\Phi^*$ one checks that the construction of the pairing in \cite[\S 4]{ta} coincides with the given in this theorem.
 Of course one has to use the $v-$ structure of $\Phi$ described  in Lemma \ref{aq} in order to define $_a\Phi^*.$
\end{proof}

\begin{rem}
In \cite{gl} the authors defined a dual t-motive in a different way cf. \cite[formula (1.8.1) for $n=1$]{gl} and showed that their definition is equivalent to that given in \cite{ta}. In \cite{hj}, in a more general context of $A$-motives, a dual of an $A$-motive  is defined by means of the internal $\hom$ for $A$-motives. According to \cite[\S 12.1]{gl3},
the aforementioned constructions  from \cite{gl} and \cite{hj} are closely related.
\end{rem}
\begin{rem}
Q. Gazda in \cite{g24} considers $\Ext^1$ functor in the category of $A$-motives. As 
in the case of a curve $P^1_{{\mathbb F}_q}$ and $A={\mathbb F}_q[t]$ abelian \tm modules correspond to effective $A$-motives, see \cite[\S 2.5]{hj}, considered by us $\Ext^1_{\tau}(\Phi,C)$ is a specialization 
of his $\Ext^1$ (cf. also \cite[\S 7]{pr}).
 Note, however, that in  Definition \ref{dual2} of a dual of a \tm module, we used $\Ext^1_{0,\tau}$ not $\Ext^1_{\tau}.$
\end{rem}
%%%%%%%%%%%%%%%%%%%%%%%%%%%%%%%%%%%%
 %%%%%%%%%%%%%%%%%%%%%%%

{}

\end{document}